\documentclass[reqno]{amsart}
\usepackage{amssymb, amscd, amsfonts}
\usepackage[mathscr]{eucal}
\usepackage{graphicx}
\usepackage{amscd}

\newtheorem{theorem}{Theorem}[section]
\newtheorem{lemma}[theorem]{Lemma}
\newtheorem{corollary}[theorem]{Corollary}
\newtheorem{proposition}[theorem]{Proposition}

\theoremstyle{definition}
\newtheorem{definition}[theorem]{Definition}

             \newtheorem{examples}[theorem]{Examples}

\newcommand{\restrict}{\,{\mathbin{\vert\mkern-0.3mu\grave{}}}\,}

\newcommand{\remove}[1]{}

\DeclareMathOperator{\McN}{\mathscr M}
\DeclareMathOperator{\McNm}{\mathscr M([0,1]^{\it m})}
\DeclareMathOperator{\McNn}{\mathscr M([0,1]^{\it n})}
\DeclareMathOperator{\conv}{\rm conv}
\DeclareMathOperator{\den}{\rm den}

\DeclareMathOperator{\id}{\rm id}

\title[Finitely
generated unital $\ell$-groups]
{Classification of finitely generated
     lattice-ordered abelian  groups with order-unit}

\author[M.Busaniche]{Manuela Busaniche $^\ddag$}
\address[M.Busaniche]{Instituto de Matem{\'a}tica Aplicada del
Litoral- FIQ, CONICET-UNL\\ Guemes 3450,
S3000GLN-Santa Fe, Argentina}
\email{manuelabusaniche@yahoo.com.ar}

\author[L.Cabrer]{Leonardo Cabrer $^\ddag$}
\address[L.Cabrer]{CONICET \\
Dep. de Matem{\'a}ticas -- Facultad de Ciencias Exactas \\
Universidad Nacional del Centro de la Provincia de Buenos Aires \\
Pinto 399 -- Tandil (7000) \\
Argentina }
\email{lcabrer@exa.unicen.edu.ar }

\author[D.Mundici]{Daniele Mundici$^\dag$}
\address[D.Mundici]{Dipartimento di
Matematica \, ``Ulisse Dini'' \\
Universit\`{a} degli Studi di Firenze \\
viale Morgagni 67/A \\
I-50134 Firenze \\
Italy}
\email{mundici@math.unifi.it }

\keywords{Lattice-ordered abelian group, rational polyhedron,
order-unit, simplicial complex, abstract simplicial complex, stellar
subdivision, Alexander starring, regular fan, De Concini-Procesi
theorem, piecewise linear function, Elliott classification, AF
$C^{*}$-algebra.}

\subjclass[2000]{Primary: 06F20, 52B11
Secondary:  52B20, 57Q15, 46L05.}

\date{\today}
\begin{document}

        \hyphenation{hom-eo-mor-phism}


\begin{abstract}
    A unital $\ell$-group $(G,u)$ is an abelian group $G$
equipped with a translation-invariant lattice-order and a
distinguished element $u$, called order-unit, whose positive integer
multiples eventually dominate each element of $G$.  We classify
finitely generated unital $\ell$-groups by sequences $\mathcal W =
(W_{0},W_{1},\ldots)$ of weighted abstract simplicial complexes, where
$W_{t+1}$ is obtained from $W_{t}$ either by the classical Alexander
binary stellar operation, or by deleting a maximal simplex of $W_{t}$.
A simple criterion is given to recognize when two such sequences
classify isomorphic unital $\ell$-groups.  Many properties of the
unital $\ell$-group $(G,u)$ can be directly read off from its
associated sequence: for instance, the properties of being totally
ordered, archimedean, finitely presented, simplicial, free.
\end{abstract}

\maketitle


\section{Introduction}
This paper deals with abelian groups $G$ equipped with a
translation-invariant lattice-order and a distinguished
{\it order-unit}
$u$, i.e., an element whose positive integer multiples eventually
dominate each element of $G$.
For brevity,  $(G,u)$
will be said to be a {\it unital $\ell$-group}.
We refer to \cite{bigkeiwol, glahol} for background.

We will classify every finitely generated unital $\ell$-group by a
sequence $\mathcal W$ of weighted abstract simplicial complexes.  A
main reason of interest in this classification is that Elliott
classification \cite{goo-af} yields a one-one correspondence $\kappa$
between isomorphism classes of unital AF $C^{*}$-algebras whose
Murray-von Neumann order of projections is a lattice, and isomorphism
classes of countable unital $\ell$-groups: $\kappa$ is an
order-theoretic enrichment of Gro\-then\-dieck $K_0$ functor,
\cite[3.9,3.12]{mun-jfa}.  Thus our Theorem
\ref{theorem:equivalence-copy} can be turned into a simple criterion
to recognize when two sequences $\mathcal W$ and $\mathcal W'$
determine isomorphic AF C*-algebras $A_{\mathcal W}$ and $A_{\mathcal
W'}$.  This criterion is reminiscent of the equivalence criterion for
Bratteli diagrams, \cite[2.7]{bra}---but our sequences are much
simpler combinatorial objects than Bratteli diagrams.

Another application of this classification immediately follows from
the categorical equivalence $\Gamma$ between unital $\ell$-groups and
MV-algebras, the algebraic counterparts of {\L}ukasiewicz
infinite-valued logic, \cite[3.9]{mun-jfa}.  Our sequences provide a
combinatorial classification of the Lindenbaum algebras of all
theories in {\L}ukasiewicz infinite-valued logic with finitely many
variables.

To describe our classifier, let us recall that a { \it (finite)
abstract simplicial complex} is a pair $H=({\mathscr V},\Sigma)\,\,$
where ${\mathscr V}$ is a finite nonempty set, whose elements are
called the {\it vertices} of $H$, and $\Sigma$ is a collection of
subsets of ${\mathscr V}$ whose union is ${\mathscr V}$, and with the
property that every subset of an element of $\Sigma$ is again an
element of $\Sigma$.  As a particular case of a construction of
Alexander \cite[p.  298]{ale}, given a two-element set $\{v,w\}\in
\Sigma$ and $a\not\in {\mathscr V}$ we define the {\it binary
subdivision} $(\{v,w\},a)$ of $H$ as the abstract simplicial complex
$(\{v,w\},a)H$ obtained by adding $a$ to the vertex set, and replacing
every set $\{v,w,u_{1},\ldots,u_{t}\}\in \Sigma$ by the two sets
$\{v,a,u_{1},\ldots,u_{t}\}$ and $\{a,w,u_{1},\ldots,u_{t}\}$ and
their faces.

A {\it weighted abstract simplicial complex} is a triple $W=({\mathscr
V},\Sigma, \omega)$ where $({\mathscr V},\Sigma)$ is an abstract
simplicial complex and $\omega$ is a map of ${\mathscr V}$ into the
set $ \{1,2,3,\ldots\}.$ For $\{v,w\} \in \Sigma$ and $a\not\in
{\mathscr V},\,$ the {\it binary subdivision} $(\{v,w\},a)W$ is the
abstract simplicial complex $(\{v,w\},a)({\mathscr V},\Sigma)$
equipped with the weight function $\tilde \omega\colon {\mathscr
V}\cup\{a\}\to \{1,2,3,\ldots\}$ given by $\tilde
\omega(a)=\omega(v)+\omega(w) \,\,\,{\rm and}\,\,\,
\tilde\omega(u)=\omega(u)$ for all $u\in {\mathscr V}.$

A sequence $\mathcal W =(W_{0},W_{1},\ldots)$ of weighted abstract
simplicial complexes is {\it stellar} if $W_{j+1}$ is obtained from
$W_{j}$ either by a deletion of a maximal set in $W_{j}$, or by a
binary subdivision, or else $W_{j+1}=W_{j}$.

We will construct a map $\mathcal W\mapsto \mathcal G(\mathcal W)$
assigning to each stellar sequence a unital $\ell$-group, and in
Theorem \ref{theorem:surject} we prove that, up to isomorphism, all
finitely generated unital $\ell$-groups arise as $\mathcal G(\mathcal
W)$ for some $\mathcal W$.  In Theorems \ref{theorem:equivalence-copy}
and \ref{theorem:equivalence-bis} we give necessary and sufficient
conditions for two stellar sequences to determine isomorphic unital
$\ell$-groups.  In Theorem \ref{theorem:tavola} we will describe the
combinatorial counterparts of such $\ell$-group properties as
se\-mi\-sim\-pli\-ci\-ty, simplicity, freeness, finite presentability,
and the property of being totally ordered.

Throughout this paper we will use techniques from polyhedral geometry,
notably regular fans and their desingularizations, \cite{ewa}.  Our
key tool is given by $\mathbb Z$-homeomorphisms (i.e.,
PL-homeomorphisms with integer coefficients).  In Theorem
\ref{theorem:leo-manu} it is proved that $\mathbb Z$-homeomorphisms
are the geometrical counterpart of unital $\ell$-isomorphisms.  All
the necessary background material will be provided in the next
section.

\section{Unital $\ell$-groups,
polyhedra and regular complexes}
\label{section:groups}
\subsection*{Lattice-ordered abelian groups with order-unit} A {\em
lattice-ordered abelian group} ({\em $\ell$-group}) is a structure
$(G,+,-,0,\vee,\wedge)$ such that $(G,+,-,0)$ is an abelian group,
$(G,\vee,\wedge)$ is a lattice, and $x+(y\vee z)=(x+y)\vee(x+z)$ for
all $x,y,z\in G$.  An {\it order-unit} in $G$ (``unit{\'e} forte'' in
\cite{bigkeiwol}) is an element $u\in G$ having the property that for
every $g\in G$ there is $0\leq n\in \mathbb Z$ such that $g\leq nu.$ A
{\it unital} $\ell$-group $(G,u)$ is an $\ell$-group $G$ with a
distinguished order-unit $u$.

The morphisms of unital $\ell$-groups, called {\it unital
$\ell$-homomorphisms}, are group-homomorphisms that also preserve the
order-unit and the lattice structure.  By an {\it $\ell$-ideal} $I$ of
$(G,u)$ we mean the kernel of a unital $\ell$-homomorphism.  Any such
$I$ determines the (quotient) unital $\ell$-homomorphism
$(G,u)\to(G,u)/I$ in the usual way \cite{bigkeiwol, glahol}.

We let $\McNn$ denote the unital $\ell$-group of piecewise linear
continuous functions $f \colon [0,1]^{n}\to \mathbb R$, such that each
piece of $f$ has integer coefficients, with the constant function 1 as
a distinguished order-unit.  The number of pieces is always finite;
``linear'' is to be understood in the affine sense.

More generally, for any nonempty subset $X\subseteq [0,1]^{n}$ we
denote by $\McN(X) $ the unital $\ell$-group of restrictions to $X$ of
the functions in $\McNn$, with the constant function $1$ as the
order-unit.

For every $\ell$-ideal $I$ of $\McNn$ we denote by $\mathcal Z(I)$ the
family of zerosets of (the functions in) $I$, in symbols,
\begin{equation}
    \label{equation:zeta} \mathcal Z(I)=\{X\subseteq
[0,1]^{n}\mid \exists g\in I \mbox{ with } X= \mathcal Z(g)=
g^{-1}(0)\}.
\end{equation}

The coordinate functions $\pi_{i}\colon [0,1]^{n} \to
\mathbb R,\,(i=1,\ldots,n)$,
together with the unit $1$,  generate
the unital  $\ell$-group $\McNn$.
They are  said to be a {\it free generating set}
of  $\McNn$ because they have the following
universal property:

\begin{theorem}
            \label{theorem:mcn}
{\rm \cite[4.16]{mun-jfa}}
Let $\{g_{1},\ldots,g_{n}\}\subseteq [0,u]$
be a set of generators of
a unital $\ell$-group $(G,u)$.
Then the {\rm map} $\pi_{i} \mapsto g_{i}$
 can be uniquely extended to a unital
$\ell$-homomorphism of $\McNn$ onto $(G,u)$.
\end{theorem}

\begin{corollary}
       \label{corollary:quotient}
Up to $\ell$-isomorphism, every finitely generated unital
$\ell$-group has the form $\McNn/I$ for some $n$ and
$\ell$-ideal $I$ of
$\McNn$.
\end{corollary}

\medskip
\subsection*{Rational polyhedra, complexes and regularity}
Following \cite[p.4]{sta}, by a {\it polyhedron} $P\subseteq \mathbb
R^{n}$ we mean a finite union of convex hulls of finite sets of points
in $\mathbb R^{n}$.  A {\it rational polyhedron } is a finite union of
convex hulls of finite sets of rational points in $\mathbb R^{n},\,
n=1,2,\ldots.$ An example of rational polyhedron $P\subseteq \mathbb
[0,1]^{n}$
is given by the zeroset $\mathcal Z(f)=f^{-1}(0)$ of any
$f\in \McNn$.  In Proposition \ref{proposition:poly} and Lemma
\ref{lemma:support} below we will see that this is the most general
possible example.


\medskip
For any rational
point $y \in \mathbb R^{n}$
we denote by   $\den(y)$
the least common denominator of the coordinates
of $y$.  The integer vector  $\tilde y =
\den(y)(y,1)\in \mathbb Z^{n+1}$
is called  the {\it homogeneous correspondent}
of $y$.
For every rational $m$-simplex
$T = \conv(v_0, \ldots, v_m) \subseteq \mathbb R^{n}$,
we will use the notation
$$
T^{\uparrow}
=\mathbb R_{\geq 0}\,\tilde v_0 +\cdots+
\mathbb R_{\geq 0}\, \tilde v_m
$$
for the positive span of $\tilde v_0,\ldots,\tilde v_m$
in $\mathbb R^{n+1}$.

We refer to \cite{ewa}
for background on simplicial complexes.
Unless otherwise specified, every
complex ${\mathcal K}$
in this paper will be simplicial,
and the adjective ``simplicial'' will be
mostly  omitted.
For every  complex  ${\mathcal K}$, its
     {\it support} $|{\mathcal K}|$ is the pointset
union of all simplexes of ${\mathcal K}$.
We say that  the complex  ${\mathcal K}$ is
{\it rational}  if
all simplexes of ${\mathcal K}$ are rational:
in this case,  the set
$$
       {\mathcal K}^{\uparrow}=
\{T^{\uparrow}\mid T\in {\mathcal K}\}
$$
is known as a simplicial fan \cite{ewa}.
A rational $m$-simplex
$T = \conv(v_0, \ldots, v_m) \subseteq \mathbb R^{n}$
is {\it regular} if the set
$\{\tilde{v}_{0},\ldots,\tilde{v}_{m}\}$
is part of a basis in the free abelian group
$\mathbb Z^{n+1}$.
A  rational complex
$\Delta$  is said to be
{\it regular}  if every simplex
$T\in \Delta$ is regular.  In other words,
the fan $\Delta^{\uparrow}$ is
regular \cite[V, \S 4]{ewa}.

We recall here some results about regular
complexes and rational
polyhedra for later use.
For detailed proofs see \cite{marmun} and
\cite{mun-dcds} where regular complexes are called ``unimodular
triangulations''.


\begin{proposition}{\rm \cite[5.1]{marmun}}
       \label{proposition:poly}
       A set   $X\subseteq [0,1]^{n}$
    coincides with the support of
	some regular complex   $\Delta$  iff
	$X=f^{-1}(0)$ for some $f\in \McNn$.
       \end{proposition}
%
%

     \begin{lemma}\cite[p.539]{mun-dcds}
	\label{lemma:support}
Any rational polyhedron
$P\subseteq [0,1]^{n}$ is the support of
some regular complex $\Delta$.
\end{lemma}


	  \remove{
\begin{proof}   A {\it rational} halfspace $H$ in
    $\mathbb R^{n}$ is a set of the form
$H =\{(x_{1},\ldots,x_{n})\in \mathbb R^{n}
     \mid  p_{1}x_{1}+\cdots + p_{n}x_{n}  \geq q\}$
     for some rational numbers
    $p_{i}$ and $q$.
    Write     $P = T_{1}\cup\ldots\cup T_{t}$
       for suitable  rational simplexes. Let
       $\mathcal H = \{H_{1},\ldots,H_{h}\}$
     be a set  of rational  half-spaces in  $\mathbb R^{n}$
such that every  $T_{j}$
       is the intersection of halfspaces
of $\mathcal H$.  Arbitrarily choose
a regular complex $\Delta_{0}$
with support $[0,1]^{n}$,
(see  \cite[p. 60-61]{sem} for one such example).
	Repeated application of the
De  Concini-Procesi
theorem  \cite{decpro},
\cite[p. 252]{ewa}, \cite{decpro},
as  is done  in
\cite[2.2]{pan-jsl},
yields a sequence  of regular complexes
$\Delta_{0},\,\, \Delta_{1},
\ldots, \Delta_{r} $  where
    each   $\Delta_{k+1}$ is obtained
by blowing-up  $\Delta_{k}$ at the Farey mediant
of some 1-simplex of  $\Delta_{k}$,  and
for each $i=1,\ldots,h$, the convex
     polyhedron  $H_{i}\cap [0,1]^{n}$
     is a union of simplexes of $\Delta_{r}.$
     It follows that each  simplex
    $T_{1},\ldots,T_{t}$
    is a union of simplexes of $ \Delta_{r} $.
    Now let $\Delta=
    \{S\in \Delta_{r}\mid S\subseteq P\}$.
       \end{proof}
       }

  \remove{
	  $(\Rightarrow)$  \
	Extend
	$\Delta$ to a regular complex $\nabla$ with support
coinciding with the	the $n$-cube, as in
{ 	Proposition \ref{proposition:extend}.} Let
$H_\nabla= \{h_v \mid v \text{ is a vertex of } \nabla \text{ not in
}\Delta\}.$
By { Proposition
	\ref{proposition:schauder-mac}},
	the required $f$ can now be obtained as
	a suitable linear combination
	of the hats in $H_{\nabla}$ with integer coefficients
	$\geq 0$.
    $\,\,(\Leftarrow)$   \
There exists a
rational
complex  $\mathcal K$ over the $n$-cube such that $f$ is linear
over each polyhedron of $\mathcal K$.
Without adding new vertices
we can refine $\mathcal K$ to a simplicial complex
$\mathcal S$, as
in \cite[III, 2.6]{ewa}.   The construction of
\cite[1.2]{mun-adv} now yields a regular complex
$\Theta$ with support  $[0,1]^{n}$ such that
    every simplex of $\mathcal S$ is
a union of simplexes of $\Theta$.
Thus $f$ is linear over each simplex of $\Theta$.
The subset $\Delta$ of $\Theta$ given by
those simplexes $S$ such that $f$ vanishes over
$S$ will provide the required regular complex with support
$X$.
}


\subsection*{Subdivision, blow-up, Farey mediant}
Given complexes ${\mathcal K}$ and ${\mathcal H}$
with $|{\mathcal K}|=|{\mathcal H}|$
we say that ${\mathcal H}$ is a
{\it subdivision} of ${\mathcal K}$
if every simplex of ${\mathcal H}$ is contained
in a simplex of ${\mathcal K}$.
For any  point $p \in |{\mathcal K}|  \subseteq \mathbb R^{n}$,
    the {\it blow-up ${\mathcal K}_{(p)}$ of ${\mathcal K}$ at
    $p$} is the    subdivision
of ${\mathcal K}$ given  by replacing
every simplex $T \in {\mathcal K}$
that contains $p$ by the set
of all  simplexes of the form $\conv(F\cup\{p\})$, where
    $F$ is any face of $T$ that does not contain $p$
     (see  \cite[III, 2.1]{ewa} or
     \cite[p. 376]{wlo}).

    For any regular $1$-simplex
    $E =\conv(v_{0},v_{1})\subseteq \mathbb R^{n}$, the {\it Farey
    mediant} of $E$ is the rational point $v$ of $E$ whose homogeneous
    correspondent $\tilde v$ coincides with $\tilde{v}_0+\tilde{v}_1$.  If
    $E$ belongs to a regular complex $\Delta$ and $v$ is the Farey mediant
    of $E$ then the blow-up $\Delta_{(v)},$ called {\it binary Farey
    blow-up}, is a regular complex.

\begin{proposition}
    \label{proposition:extended-dcp}
    Suppose we are
given rational polyhedra $Q\subseteq P\subseteq [0,1]^{n}$ and a
regular complex $\Delta$ with support $P$.  Then there is a
subdivision $\Delta^{\natural}$ of $\Delta$ obtained by binary Farey
blow-ups, such that $Q=\bigcup\{T\in \Delta^{\natural}\mid T\subseteq
Q\}.$
\end{proposition}

     \begin{proof} We closely follow the argument of the proof in
     \cite[p.539]{mun-dcds}.  Let us write $Q = T_{1}\cup\ldots\cup T_{t}$
     for suitable rational simplexes.  Let $\mathcal H =
     \{H_{1},\ldots,H_{h}\}$ be a set of rational half-spaces in $\mathbb
     R^{n}$ such that every $T_{j}$ is the intersection of halfspaces of
     $\mathcal H$.  Using the De Concini-Procesi theorem \cite[p.252]{ewa},
     we obtain a sequence of regular complexes $\Delta=\Delta_{0},\,\,
     \Delta_{1}, \ldots, \Delta_{r} $ where each $\Delta_{k+1}$ is obtained
     by blowing-up $\Delta_{k}$ at the Farey mediant of some 1-simplex of
     $\Delta_{k}$, and for each $i=1,\ldots,h$, the convex polyhedron
     $H_{i}\cap [0,1]^{n}$ is a union of simplexes of $\Delta_{r}.$ It
     follows that each simplex $T_{1},\ldots,T_{t}$ is a union of simplexes
     of $ \Delta_{r} $.  Now $\Delta^{\natural}=\Delta_{r}$ yields the
     desired subdivision of $\Delta$.  \end{proof}


       The following proposition states
       that every $\ell$-ideal $I$ of $\McN(P)$ is uniquely
       determined by the zerosets of all functions in $I$:

        \begin{proposition}
\label{proposition:woj}
Let $P\subseteq [0,1]^{n}$ be a rational polyhedron
and  $I$  an $\ell$-ideal  of $\McN(P)$.
Then   for every
$f\in \McN(P)$ we have:  $f\in I$ iff
$f^{-1}(0)\supseteq g^{-1}(0)$ for some $g\in I$.
\end{proposition}
\begin{proof} For the nontrivial direction,
suppose  $f^{-1}(0)\supseteq g^{-1}(0)$
       and without loss of generality,
       $g \geq 0,$  and $f\geq 0.$
       We must find $0\leq m\in \mathbb Z$
    such that  $mg \geq f$.
    Since $f,g\in \McN(P)$, by suitably
    subdividing the regular complex
    given by { Lemma \ref{lemma:support}}
    we obtain a 
    complex  $\Delta$   all of whose simplexes
    $T_{1},\ldots,T_{s}$
    have rational vertices,  with $\bigcup T_{i}
    = P$  and such that over every
    $T_{i}$  both $f$ and $g$ are linear. Fix $T_{i}\in\Delta,$ since
    for every vertex $v$ of $T_{i}$,
    $f(v)\neq 0$ implies $g(v)\neq0$, it is easy
    to find and integer $m_{i}>0$ such that
    $m_{i}g \geq f$ over  $T_{i}$.
    Now let $m=\max(m_{1},\ldots,m_{u}).$
       \end{proof}

Further immediate
properties of zerosets are recorded
here for later use:

\begin{proposition}
       \label{proposition:zeroset}
       If $I$ is an $\ell$-ideal of $\McNm,\,\,$
       $P\in \mathcal Z(I)$ and $0\leq f$
       is the restriction to $P$  of some $k\in I$ (in symbols
       $f=k\restrict P$),
       then there is an extension
$\tilde f$ of $f$  such that
$\tilde f \in I$  and
$ \tilde f>0$ over $[0,1]^{m}\setminus P$.
Thus $\mathcal Z(f)=\mathcal Z(\tilde f)$
for some $\tilde f\in I$.
\end{proposition}


\begin{proposition}
       \label{proposition:threesets}
Let  $I$ be an $\ell$-ideal of $\McNm$
and $P\in \mathcal Z(I)$.  Define
\begin{itemize}
\item[(i)]
       $I\restrict P=
\{f\restrict P\mid f\in I\},$

\smallskip
\item[(ii)]
$\mathcal Z(I)_{\cap P}  =  \{X\cap P\mid X\in \mathcal Z(I)\},$

\smallskip
\item[(iii)]
$\mathcal Z(I)_{\subseteq P}=\{X\in
\mathcal Z(I)\mid X\subseteq P\}$.
\end{itemize}
Then $\mathcal Z(I\restrict P)=\mathcal Z(I)_{\cap P}
=\mathcal Z(I)_{\subseteq P}$.
    \end{proposition}


    %


       \section{$\mathbb Z$-hom\-e\-o\-mor\-phism
       and unital
       $\ell$-isomorphism}
       \label{section:zeta}
            Given  rational polyhedra
       $P\subseteq \mathbb{R}^m$ and $Q \subseteq
\mathbb{R}^n$, a piecewise linear homeomorphism
$\eta$ of $P$ onto $Q$ is
said to be
a {\it $\mathbb Z$-hom\-e\-o\-mor\-phism},
in symbols, $\eta\colon
P\cong_{\mathbb Z} Q,$ if all linear pieces
of $\eta$ and $\eta^{-1}$
have integer coefficients.

       \medskip
The following first main result
of this paper highlights
the mutual relations between
       $\mathbb Z$-hom\-e\-o\-mor\-phisms of polyhedra and
       unital  $\ell$-isomorphisms of finitely generated
       unital $\ell$-groups,
       as represented by { Corollary
       \ref{corollary:quotient}}:

         \begin{theorem}
        \label{theorem:leo-manu}
For any  $\ell$-ideals $I$ of
        $\McNm$ and $J$ of $\McNn$
        the following conditions are equivalent:

        \begin{itemize}
\item[(i)]
    $ \McNm/I
\cong \McNn/J$.

\smallskip
\item[(ii)]
    For some $P\in \mathcal Z(I)$,
     $Q\in \mathcal Z(J)$ and ${\mathbb Z}$-homeomorphism
     $\eta$ of $P$ onto $Q$,
     the map $X\mapsto \eta(X)$
     sends $\mathcal Z(I)_{\cap P} $
     one-one onto
      $\mathcal Z(J)_{\cap Q}.$

     \end{itemize}
        \end{theorem}

        \begin{proof}
$(i)\to(ii)$  For definiteness, let us write
$\iota\colon
\McNm/I\cong \McNn/J,$
and  $\epsilon = \iota^{-1}.$
%
Let $\id_{m}$ denote the identity $(\pi_{1},\ldots,\pi_{m})$ over the
$m$-cube, and $\id_{n}$ the identity over the $n$-cube.  Each element
$\pi_{i}/I\in \McNm/I$ is sent by $\iota$ to some element $a_{i}/J$ of
$\McNn/J.$ Writing $[0,1] \ni ((a_{i}/J)\vee 0)\wedge 1 =((a_{i}\vee
0)\wedge 1)/J, $ and replacing, if necessary, $a_{i}$ by $(a_{i}\vee
0)\wedge 1$, it is no loss of generality to assume that $a_{i}$
belongs to the unit interval of $\McNn,$ i.e., the range of $a_{i}$ is
contained in the unit interval $[0,1]$.  Thus for a suitable $m$-tuple
$a=(a_{1},\ldots,a_{m})$ of functions $a_{i}\in \McNn$ we have $a
\colon [0,1]^{n} \to [0,1]^{m}$.  Symmetrically, for some
$b=(b_{1},\ldots,b_{n}) \colon [0,1]^{m} \to [0,1]^{n}, \,\,\,
b_{j}\in \McNm$ we can write
	     \begin{equation}
	     \label{equation:two-bis}
\iota\colon \id_{m}/I\mapsto a/J\,\,\,
\mbox{and}\,\,\,
\epsilon \colon \id_{n}/J\mapsto b/I.
	     \end{equation}
%
%
%
For any   $f \in \McNm$
and $g\in \McNn,$
arguing  by induction on the
number of operations in
$f$ and $g$ in the light of
{ Theorem \ref{theorem:mcn}},
we get the following
generalization of (\ref{equation:two-bis}):

	     \begin{equation}
	     \label{equation:five-bis}
	     \iota\colon f/I\mapsto (f\circ a)/J
	     \,\,\,\mbox{and}\,\,\,
\epsilon \colon g/J\mapsto (g\circ b)/I.
	     \end{equation}
%
%
    %
%
%
It follows that
$$
\frac{f}{I}= (\epsilon\circ\iota)\frac{f}{I}=
\epsilon(\iota(\frac{f}{I}))
= \epsilon(\frac{f\circ a}{J}) = \frac{f\circ a\circ b}{I}.
$$
By definition of the congruence
induced by $I$, the  function
$|\pi_i-a_i\circ b|=|\pi_i-\pi_i\circ a\circ b|$
belongs to $I,\,\,\,(i=1,\ldots  ,m)$.
Here, as usual,  $|\cdot|$ denotes
absolute value.
%
%
It follows that the function $e = \sum_{i=1}^{m}\,\,|\pi_{i} -\,
a_{i}\circ b\,|$ belongs to $I$, and its zeroset $\mathcal Z(e)$
belongs to $\mathcal Z(I)$.  The set $P=\mathcal Z(e)$ satisfies the
identity $$P=\{x\in [0,1]^{m}\mid (a\circ b)(x)=x\}.$$ One similarly
notes that the set $Q=\{y\in [0,1]^{n}\mid (b\circ a)(y)=y\}$ belongs
to $\mathcal Z(J).$ By construction, the restriction of $b$ to $P$
provides a $\mathbb Z$-hom\-eo\-mor\-phism $\eta$ of $P$ onto $Q$,
whose inverse $\theta$ is the restriction of $a$ to $Q$.  In symbols,
$$ b\restrict P=\eta\colon P\cong_{\mathbb Z} Q, \quad a \restrict
Q=\theta\colon Q\cong_{\mathbb Z} P. $$
     Suppose $Y\in  \mathcal Z(J)_{\cap Q}$,
     with the intent of proving that
     $\theta(Y)\in \mathcal Z(I)_{\cap P}.$
By { Proposition
\ref{proposition:zeroset}},
we can  write  $Y = \mathcal Z(k) =
\mathcal Z(k\restrict Q)$ for some
$0 \leq  k\in J.$
By  (\ref{equation:five-bis}),
    the composite function $k\circ b$ belongs to $I.$
%
Since $e \in I$, $k \circ b + e\in I$.  Because $\mathcal Z(e)=P$, the
function $k \circ b + e$ coincides with $k \circ b$ over $P$, and we
have $$ \theta(Y) =\theta(\mathcal Z(k\restrict Q))= \mathcal
Z((k\restrict Q) \circ \eta) =\mathcal Z(k \circ \eta)=\mathcal Z(k
\circ b \restrict P)= $$ $$ \mathcal Z(k \circ b + e)= P \cap \mathcal
Z(k\circ b)\in \mathcal{Z}(I)_{\cap P} $$ as desired.  Reversing the
roles of $\theta$ and $\eta$, we have the required one-one
correspondence $X\mapsto \eta(X)$ of $\mathcal Z(I)_{\cap P} $ onto
$\mathcal Z(J)_{\cap Q}$.

    \bigskip
        $(ii)\to(i)\,\,\,$
Let  $I_{P} $ (resp., let $J_{Q}$) be the $\ell$-ideal of
$\McNm$ (resp., of $\McNn$) given by those functions that
vanish over $P$ (resp., over $Q$).
By \cite[5.2]{marmun} we have
unital
$\ell$-isomorphisms
\begin{equation}
       \label{equation:wellknown}
       \alpha\colon \McN(P) \cong \McNm/I_{P}
    \quad \mbox{ with }\quad  \alpha(I\restrict P)=I/I_{P}
\end{equation}
and
\begin{equation}
       \label{equation:wellknown-bis}
       \beta\colon \McN(Q)\cong\McNn/J_{Q}
         \quad \mbox{ with }\quad  \beta(J\restrict Q)=J/J_{Q}.
       \end{equation}
As a particular case of a general algebraic result,
the map
$$
\frac{f/I_{P}}{I/I_{P}}\mapsto \frac{f}{I}
$$
is a unital $\ell$-isomorphism
of  $\frac{\McNm/I_{P}}{I/I_{P}}$ onto
$\frac{\McNm}{I}.$
{}From  (\ref{equation:wellknown})-(\ref{equation:wellknown-bis})
we have
unital $\ell$-isomorphisms
\begin{equation}
\label{equation:uno}
\frac{\McNm}{I} \cong \frac{\McNm/I_{P}}{I/I_{P}}\cong
\frac{\McN(P)}{I\restrict P}
\end{equation}
and
\begin{equation}
\label{equation:uno-bis}
\frac{\McNn}{J} \cong \frac{\McNn/J_{Q}}{J/J_{Q}}\cong
\frac{\McN(Q)}{J\restrict Q}.
\end{equation}
{}Letting $\theta=\eta^{-1},$  we have
$
\theta \colon Q\cong_{\mathbb Z}  P
$
and the map $ \lambda\colon
k\mapsto k\circ\theta$ is a
    unital $\ell$-isomorphism
of $\McN(P)$  onto
$\McN(Q)$.
Further, the map  $Y\mapsto \theta(Y)$  sends
$
\mathcal Z(J)_{\cap Q}
=\mathcal Z(J\restrict Q)$
one-one onto
$
\mathcal Z(I)_{\cap P}
=\mathcal Z(I\restrict P).
$

\bigskip
\noindent
{\it Claim.} The restriction of $\lambda$ to the $\ell$-ideal $
I\restrict P$ of $\McN(P)$ maps $ I\restrict P$ one-one onto
$J\restrict Q$.  Thus the map $$ \frac{k}{I\restrict P} \mapsto
\frac{\lambda(k)}{\lambda(I\restrict P)} $$ defines a unital
$\ell$-isomorphism of $\McN(P)/(I\restrict P)$ onto
$\McN(Q)/(J\restrict Q)$.

\smallskip
By Proposition \ref{proposition:threesets},
for each $l\in \McN(P)$  if
$l\in I\restrict P$ then $\mathcal Z(l)
\in \mathcal Z(I\restrict P)
=\mathcal Z(I)_{\cap P}$. Thus by definition of $\lambda$,
$\mathcal Z(\lambda(l))
=\mathcal Z(l \circ\theta)
= \eta(\mathcal Z(l))\in \mathcal Z(J)_{\cap Q}$, i.e.,
$\lambda(l)\in J\restrict Q.$

Reversing
    the roles of $\lambda$ and $\lambda^{-1}$,
    our claim is settled.

\bigskip
From (\ref{equation:uno}) it now follows that the map $$
\frac{g}{I}\mapsto \frac{g\restrict P}{I\restrict P} \mapsto
\frac{\lambda(g\restrict P)}{\lambda(I\restrict P)} =
\frac{(g\restrict P)\circ\theta}{J\restrict Q} $$ is a unital
$\ell$-isomorphism of $ {\McNm}/{I}$ onto $ {\McN(Q)}/({J\restrict
Q})$.  {}From (\ref{equation:uno-bis}) we obtain the desired
conclusion.  \end{proof}


We will make use of the following

\begin{lemma}
\label{lemma:dominio_lineal}
Let $\Delta $ be a regular complex whose support is contained in the
$n$-cube. Let  $\eta$ be a $\mathbb{Z}$-homeomorphism defined over
$|\Delta|.$ Then there is a regular complex $\Delta'$ obtained from
$\Delta$ by binary Farey blow-ups, such that $\eta$ is linear over each
simplex of $\Delta'.$
\end{lemma}

\begin{proof} Since $\eta$ has only finitely many pieces
      and each piece is linear with integer coefficients,
      a routine argument yields a  (rational, simplicial)
      complex $\nabla$ with support  $|\Delta|$ such that
      $\eta$  is linear over every simplex of $\nabla$.
Now iterated application of { Proposition
\ref{proposition:extended-dcp}} yields the desired
regular complex $\Delta'$.
\end{proof}

For later purposes,  we record here the following
trivial
property of {\it linear}
   $\mathbb  Z$-homeo\-mor\-phisms:

\begin{lemma}
       \label{lemma:un-simplesso-jem}
       Let
$T=\conv(v_{0},\ldots,v_{k})\subseteq \mathbb R^{m}$ and
$U=\conv(w_{0},\ldots,w_{k})\subseteq \mathbb R^{n}$ be regular
$k$-simplexes.
If  $\den(v_{i})=\den(w_{i})$ for all
$i=0,\ldots,k$,  then there is precisely one    linear
$\mathbb Z$-homeomorphism
$\eta_{T}$ of $T$ onto $U$ such that
$\eta_{T}(v_{i})=w_{i}$ for all $i$.
\end{lemma}

\begin{lemma}
         \label{lemma:image_homeo}
         Let $S\subseteq \mathbb{R}^n$ be a   regular simplex
and $\eta \colon S \rightarrow\mathbb{R}^m$ a
$\mathbb{Z}$-hom\-eo\-morph\-ism which is linear over $S.$
   Let $\eta(S)=S'.$
   Then $S'$
is a   regular simplex.
\end{lemma}


\begin{proof}  For some integer vector
 $b=(b_1,\ldots,b_m)\in \mathbb Z^{m}$ and  $m\times n$
 integer matrix  $\mathcal M$ we have
 $\eta(x)=\mathcal  M x+b$ for each $x\in S$.
 Let  $r_{i1},\ldots,r_{in}$ be the $i$th row of $\mathcal  M$.
Let   $\mathcal{M}_S$ be the
 $(m+1)\times (n+1)$ integer matrix whose
 $i$th  row is $(r_{i1},\ldots,r_{in},b_i)$,
for $i=1, \ldots, m$,   and whose
 bottom row
 is  ${(0,\ldots, 0,1)}.$
The homogeneous
linear map $(x,1) \mapsto \mathcal{M}_S (x,1)$ sends each vector
$(x,1)\in S\times\{1\}$ into a vector $(y,1)\in S'\times\{1\}.$ Let
$\tilde{v}_0, \ldots , \tilde{v}_j\in \mathbb{Z}^{n+1}$ be the
homogeneous correspondents of the vertices $v_0, \ldots , v_j$ of $S$
and let
$$ S^{\uparrow} =\mathbb R_{\geq 0}\,\tilde v_0 +\cdots+ \mathbb
R_{\geq 0}\, \tilde v_j $$ be the positive span of $\tilde{v}_0,
\ldots , \tilde{v}_j.$ Similarly let $S'^{\uparrow}$ be the positive
span in $\mathbb{R}^{m+1}$ of the integer vectors $\mathcal{M}_S \
\tilde v_0, \ldots , \mathcal{M}_S \ \tilde{v}_j.$ By construction,
$\mathcal{M}_S$ sends integer points of $S^{\uparrow}$ one-one into
integer points of $S'^{\uparrow}.$ Interchanging the roles of $\eta$
and $\eta^{-1}$, we see that $\mathcal{M}_S$ actually sends the set of
integer points in $S^{\uparrow}$ one-one {\it onto} the set of
integers points of $S'^{\uparrow}.$ By Blichfeldt theorem,
\cite{grulek}, the regularity of $S$ is equivalent to saying that the
half-open parallelepiped $Q_S=\{ \mu_0\tilde{v}_0+ \ldots
\mu_j\tilde{v}_j |0 \le \mu_0, \ldots \mu_j <1 \}$ contains no nonzero
integer points.  One also has a similar characterization of the
regularity of $S'$.  The mentioned properties of $\mathcal{M}_S,$
together with the assumed regularity of $S$ ensure that $S'$ is
regular.
\end{proof}

   \begin{corollary} \label{corollary:reg_into_reg} Let $\Delta$ be a
   regular complex and $\eta$ a $\mathbb{Z}$-homeomorphism of $|\Delta|$
   which is linear over each simplex of $\Delta.$ Then the set
   $\eta(\Delta)=\{\eta(S) \mid S\in \Delta\}$ is a regular complex.
   \end{corollary}

\section{Realizations of weighted abstract simplicial complexes}
We now turn to the combinatorial  objects
defined in the Introduction.
For every regular complex $\Lambda,$
   the {\it skeleton} of  $\Lambda$ is
the  weighted abstract simplicial complex
$W_{\Lambda}=(\mathcal {\mathscr V},\Sigma,\omega)$
given by the following stipulations:
\begin{enumerate}
      \item  $\mathcal {\mathscr V}=$ vertices of $\Lambda$.
      \item   For every vertex $v$ of $\Lambda$,
$\omega(v)=\den(v).$
\item For every subset
$W=\{w_{1},\ldots,w_{k}\}$
of ${\mathscr V}$,
$W\in \Sigma$ iff $\conv(w_{1},\ldots,w_{k})\in \Lambda.$
\end{enumerate}

Given two weighted abstract simplicial
complexes $ W= (\mathcal V,\Sigma,\omega)$ and
$W' = (\mathcal V',\Sigma',\omega') $ we
write
$$\gamma\colon W \cong W',$$
    and we say that $\gamma$ is a
     {\it  combinatorial isomorphism} between
     $W$  and $W'$,
if $\gamma$  is a one-one map
    from $\mathscr V$ onto
    ${\mathscr V}'$ such that $\omega'(\gamma(v))
    =\omega(v)$  for all  $v\in {\mathscr V}$, and
    $\{w_{1},\ldots,w_{k}\}\in \Sigma$
    iff $\{\gamma(w_{1}),\ldots,\gamma(w_{k})\}\in \Sigma'$
    for each subset  $\{w_{1},\ldots,w_{k}\}$ of $\mathscr V$.

\begin{definition} \label{definition:delta-realization} Let $W$ be a
weighted abstract simplicial complex and $\nabla$ a regular complex.
Then a {\it $\nabla$-realization} of $W$ is a combinatorial
isomorphism $\iota$ between $W$ and the skeleton $W_{\nabla}$ of
$\nabla$.  We write $\iota\colon W\to \nabla$ to mean that $\iota$ is
a $\nabla$-realization of $W$.  \end{definition}

\subsection*{Examples of realizations}
For any regular complex
$\Lambda$,  the
    identity function over the set of vertices of
    $\Lambda$  is
a  $\Lambda$-realization of
$W_{\Lambda},$  called the {\it trivial
realization of the skeleton}  $W_{\Lambda}$.

Symmetrically, let $W=({\mathscr V},\Sigma, \omega)$
be  a  weighted abstract simplicial complex
with vertex set  ${\mathscr V}=\{v_{1},\ldots,v_{n}\}.$
For   $e_{1},\ldots,e_{n}$ the standard basis vectors
of $\mathbb R^{n}$,  let
    $\Delta_{W}$ be the complex  whose
vertices are
$$v'_{1}=
e_{1}/\omega(v_{1}),\ldots,v'_{n}=e_{n}/\omega(v_{n}),
$$
and whose $k$-simplexes
($k=0,\ldots,n$) are given by
$$
\conv(v'_{i(0)},\ldots,
v'_{i(k)})\in \Delta_{W} \text{ iff }
\{ v_{i(0)},\ldots, v_{i(k)}\}\in \Sigma.
$$
Note that $\Delta_{W}$ is a regular
complex  and
$|\Delta_{W}|\subseteq [0,1]^{n}$.
The function
\begin{equation}
\label{equation:iota}
\tilde{\iota}\colon
v_{i}\in \mathcal V\mapsto v'_{i}
\in [0,1]^{n}
\end{equation}
is
a $\Delta_{W}$-realization of $W$, called the
{\it canonical  realization} of $W$.
The dependence on the order in which
the elements $\{v_{1},\ldots,v_{n}\}$
are listed,  is tacitly understood.

\smallskip
Combining
      { Lemma \ref{lemma:un-simplesso-jem}} with a routine
    patching argument we easily get

\begin{lemma}
       \label{lemma:unione-jems}
       Let $\Lambda$ and $\nabla$ be
regular complexes, with $|\Lambda|\subseteq \mathbb R^{m}$ and
$|\nabla|\subseteq \mathbb R^{n}$.  We then have:

\begin{itemize}
\item[(i)] If  $\theta\colon W_{\Lambda}\cong W_{\nabla}$ is a
combinatorial isomorphism between the skeletons
of  $\Lambda$ and $\nabla$,  then
there is a $\mathbb Z$-homeomorphism
$\eta_{\theta}$ of $|\Lambda|$ onto $|\nabla|$ such that
$\eta_{\theta}(v)=\theta(v)$
for each vertex $v$ of $\Lambda$,  and $\eta_{\theta}$  is linear
over each simplex of $\Lambda$.

\smallskip
\item[(ii)] Letting, in particular,
$\nabla = \Delta_{W_{\Lambda}}$,
it follows that the
combinatorial isomorphism $\tilde{\iota}$ of
(\ref{equation:iota}) between
    $W_{\Lambda}$ and $W_{\nabla}$
uniquely extends to a
${\mathbb  Z}$-homeomorphism  $\eta_{\tilde{\iota}}$  of
$|\Lambda|$ onto
$|\nabla|$
such that  $\eta_{\tilde{\iota}}$
is linear over each simplex of $\Lambda$.
\end{itemize}
\end{lemma}

\medskip
\subsection*{Filtering and confluent sets of rational polyhedra}
Given a set $E$, a non\-empty family $\mathcal P$ of subsets of $E$ is
said to be {\it filtering} if whenever $X,Y\in \mathcal P$ there is
$Z\in \mathcal P \mbox{ such that } Z\subseteq X\cap Y.$ If $\mathcal
Q$ is another filtering sequence of subsets of $E$, we say that
$\mathcal P$ {\it minorizes} $\mathcal Q$, if for all $A\in \mathcal
Q$ there is $B\in \mathcal P \mbox{ such that } B\subseteq A.$ If
$\mathcal P$ and $\mathcal Q$ happen to minorize each other, we say
that $\mathcal P$ is {\it confluent} with $\mathcal Q$, and we write
$$\mathcal P\uparrow \mathcal Q.$$

\begin{lemma} \label{lemma:ideal} If $\mathcal P$ is a filtering
family of rational polyhedra in the $n$-cube then the set $$
\langle\mathcal P\rangle = \{f\in \McNn \mid \mathcal{Z}(f)\supseteq Q
\mbox{ for some } Q \in \mathcal P\} $$ is an $\ell$-ideal of $\McNn$.
\end{lemma}

\begin{proof} The zero function belongs to
       $\langle\mathcal P\rangle$ because
       $\mathcal P$  is nonempty.
       If  $0\leq f\in \langle\mathcal P\rangle$
       and $0\leq g \leq f$ then
       $g\in \langle\mathcal P\rangle$ because
       $\mathcal Z(g)\supseteq \mathcal Z(f).$ Finally, from
       $0\leq f,  0\leq g$  and $f,g\in \mathcal P$
       it follows that $f+g\in \mathcal P$ because
       $\mathcal Z(f+g)=\mathcal  Z(f)
       \cap \mathcal  Z(g)\supseteq X\cap Y$
       for some $X,Y\in \mathcal P$. Thus
       $\mathcal  Z(f+g)\supseteq Q$
       for some $Q\in \mathcal P$.
       \end{proof}

Confluence  is an equivalence
relation between  filtering
families $\mathcal P$  of  rational polyhedra in
the $n$-cube.  The equivalence
classes are in one-one correspondence with
$\ell$-ideals of $\McNn$  via the map
$\mathcal P\mapsto \langle\mathcal P\rangle.$
Specifically, we have the following
result, which will be used without explicit mention
throughout the rest of this paper:

	\begin{proposition}
\label{proposition:tre}
Given filtering families $\mathcal
P$  and  $\mathcal Q$ of rational
polyhedra in the same $n$-cube, the following
conditions are equivalent:
\begin{itemize}
       \item[(i)] $\langle\mathcal P\rangle= \langle\mathcal
Q\rangle$.
\item[(ii)] $\mathcal P\uparrow \mathcal Q$.
\item[(iii)]
$\mathcal \mathcal \mathcal Z(\langle\mathcal P\rangle) =\mathcal
\mathcal \mathcal Z(\langle\mathcal Q\rangle)$.
\end{itemize}
\end{proposition}

\begin{proof}
       $(iii)\leftrightarrow (i)$
       By { Proposition \ref{proposition:woj}},
       the $\ell$-ideals  $\langle\mathcal
P\rangle$ and $\langle\mathcal Q\rangle$
are uniquely determined by their zerosets.
$(i)\rightarrow (ii)$ Suppose $A\in \mathcal P$ is such that
there is no $B\in \mathcal Q$ such that $B\subseteq A.$ By
{ Proposition
\ref{proposition:poly} and Lemma \ref{lemma:support}}, there is $f\in
\McNn$ with $\mathcal Z(f)=A,$ and hence, $f\in \langle\mathcal
P\rangle.$ Our hypothesis about $A$ ensures that $f \not\in
\langle\mathcal Q\rangle$.
$(ii)\rightarrow (i)$ is
trivial.
\end{proof}


\subsection*{Stellar transformations}
Let $W=(\mathcal V, \Sigma,\omega)$
and $W'$ be two weighted abstract simplicial complexes. A map
$\flat: W\rightarrow W'$ is called a
{\it stellar transformation} if
$\flat$ is either a deletion of a maximal set of $\Sigma$,
or a binary subdivision,
(as defined in the Introduction)
or else $\flat$  is  the identity map.

\medskip
Recalling
Definition \ref{definition:delta-realization}
we have

\begin{lemma}
\label{lemma:meno-star}
Let $W=(\mathcal V, \Sigma,\omega)$
and $W'=(\mathcal V ', \Sigma ',\omega ')$ be two weighted abstract
simplicial complexes,
$\Delta$ a regular complex, and $\iota$  a
$\Delta$-realization of $W$,
$\iota\colon W \to \Delta$.
    Suppose that $\flat\colon W\to W'$ is a stellar transformation.
\begin{itemize}
       \item[(i)]
       In case $\flat$  deletes a maximal set
     $M\in \Sigma$, let
       $\flat ({\iota})\colon \Delta\to\Delta'$
       delete from $\Delta$ the
corresponding maximal simplex $\conv(\iota(M))$.
Then the map
$\iota'=\iota\restrict \mathcal V '$
is a $\Delta'$-realization of
$W'$.

\smallskip
\item[(ii)]
In case $\flat$  is the binary subdivision $W'=(\{a,b\}c)W$
at some two-element set $E=\{a,b\}\in \Sigma,$
and $c\not\in \mathcal V$, let
$e$ be the Farey mediant of the $1$-simplex
$\conv(\iota(E))$. Let $\flat({\iota})$  be the
Farey blow-up
    $\Delta'=\Delta_{(e)}$ of $\Delta$ at $e$.
Then  the map $\iota'=\iota\cup\{(c,e)\}$
is a $\Delta'$-realization of $W'.$
   \end{itemize}

   \medskip
   \noindent  Further, we have a
    commutative diagram
$$
\def\normalbaselines{\baselineskip20pt
\lineskip3pt
\lineskiplimit3pt}
\def\mapright#1{\smash{\mathop{\longrightarrow}\limits^{#1}}}
\def\mapdown#1{\Big\downarrow\rlap{$\vcenter{\hbox{$\scriptstyle#1$}}$}}
\begin{matrix}
W &\mapright{\flat}&W'\cr
\mapdown{\iota}& &\mapdown{\iota'}\cr
     \Delta&\mapright{\flat({\iota})}&  \Delta'\cr
\end{matrix}
$$

\noindent We say that
$\flat(\iota)$ is  {\rm the $\Delta$-transformation of} $\flat$.
(It is tacitly understood that
if  $\flat$ is the identity map, then  $\flat(\iota): \Delta \to
\Delta'$ is the identity function.)
\end{lemma}


\section{Classification of unital $\ell$-groups}
In this section we will construct
a map $\mathcal W\mapsto \mathcal G(\mathcal W)$,
from stellar sequences
to unital $\ell$-groups and prove that
the map  is onto all
finitely generated unital $\ell$-groups.

\subsection*{Main Construction}
Let $\mathcal W=W_{0},W_{1},\ldots$ be a stellar sequence.  For each $
j=0,1,\ldots\,$ let $\flat_{j}$ be the corresponding stellar
transformation sending $W_{j}$ to $W_{j+1}.$ For some $n\geq 1$ and
regular complex $\Delta_{0}$ in the $n$-cube let $\iota_{0}$ be a
$\Delta_{0}$-realization of $W_{0}.$ Then { Lemma
\ref{lemma:meno-star}} yields a commutative diagram

     \begin{equation}
       \label{equation:infinito}
\def\normalbaselines{\baselineskip20pt
\lineskip3pt
\lineskiplimit3pt}
\def\mapright#1{\smash{\mathop{\longrightarrow}\limits^{#1}}}
\def\mapdown#1{\Big\downarrow
\rlap{$\vcenter{\hbox{$\scriptstyle#1$}}$}}
\begin{matrix}
W_{0}  & \mapright{\flat_{0}} &
W_{1}  & \mapright{\flat_{1}} &
W_{2}  & \ldots\cr
\mapdown{\iota_{0}} &  &
\mapdown{\iota_{1}} &  &
\mapdown{\iota_{2}} & \cr
     \Delta_{0} & \mapright{\flat_0(\iota_{0})} &
     \Delta_{1} & \mapright{\flat_1(\iota_{1})} &
     \Delta_{2} & \ldots\cr
\end{matrix}
\end{equation}

%
%
%
%
%
%

The sequence of supports $|\Delta_{0}| \supseteq |\Delta_{1}|
\supseteq \cdots$ is called the {\it $\Delta_{0}$-orbit of $\mathcal
W$} and is denoted $\mathcal O(\mathcal W,\Delta_{0})$ (the role of
$\iota_{0}$ being tacitly understood).  As in { Lemma
\ref{lemma:ideal},} the filtering set $\mathcal O(\mathcal
W,\Delta_{0})$ determines the $\ell$-ideal $\mathcal I(\mathcal
W,\Delta_{0}) = \langle\mathcal O(\mathcal W,\Delta_{0})\rangle$ of
$\McNn$, as well as the unital $\ell$-group $ \mathcal G(\mathcal
W,\Delta_{0}) = \McNn/\mathcal I(\mathcal W,\Delta_{0}).  $ In the
particular case when $\iota_{0}$ is the canonical realization of
$W_{0}$ we write $$\mathcal O(\mathcal W),\,\,\, \mathcal I(\mathcal
W),\,\,\, \mathcal G(\mathcal W)$$ instead of $\mathcal O(\mathcal W,
\Delta_{W_{0}}),\,\,\, \mathcal I(\mathcal W, \Delta_{W_{0}}),\,\,\,
\mathcal G(\mathcal W, \Delta_{W_{0}}).  $

\smallskip
\noindent Various examples of the map $\mathcal W\mapsto
\mathcal G(\mathcal W)$ will be given in
    \ref{examples:table}  below.

    \smallskip
Our second main result is the following

\begin{theorem}
       \label{theorem:surject}
       For every
finitely generated unital $\ell$-group $(G,u)$ there is a stellar
sequence $\mathcal W$ such that $\mathcal G(\mathcal W) \cong (G,u)$.
\end{theorem}

As a preliminary step for the proof
we need the following immediate consequence
of the definitions:

\begin{lemma}
\label{lemma:tre-sei}
For any weighted abstract simplicial complex
$W$ and regular complexes
$\nabla$ and $\Delta$, let
$\iota$  be
a $\nabla$-realization of $W$,
and   $\epsilon$
a $\Delta$-realization of  $W$.
Let $\eta_{\gamma}\colon |\nabla|\to
|\Delta|$
be the $\mathbb{Z}$-homeomorphism of { Lemma
\ref{lemma:unione-jems}} corresponding to the combinatorial
isomorphism $\gamma=\epsilon\circ\iota^{-1}$.
Suppose the stellar transformation $\flat$ transforms $W$
into $W'$.  Let the commutative diagram

\begin{equation}
      \label{equation:restored}
\def\normalbaselines{\baselineskip20pt
\lineskip3pt
\lineskiplimit3pt}
\def\mapright#1{\smash{\mathop{\longrightarrow}\limits^{#1}}}
\def\mapdown#1{\Big
\downarrow\rlap{$\vcenter{\hbox{$\scriptstyle#1$}}$}}
\def\mapup#1{\Big
\uparrow\rlap{$\vcenter{\hbox{$\scriptstyle#1$}}$}}
\begin{matrix}
      \Delta &\mapright{\flat(\epsilon)}&\Delta'\cr
\mapup{\epsilon}& &\mapup{\epsilon'}\cr
W  &\mapright{\flat}& W'\cr
\mapdown{\iota}& &\mapdown{\iota'}\cr
\nabla &\mapright{\flat(\iota)}&\nabla'\cr
\end{matrix}
\end{equation}
be as given by  {   Lemma
\ref{lemma:meno-star}(ii)}.  Let further
$\gamma'=\epsilon'\circ\iota'^{-1}$,
and $\eta_{\gamma'}$ be the $\mathbb{Z}$-homeomorphism
of $|\nabla'|$ onto $|\Delta'|$ given by
{Lemma \ref{lemma:unione-jems}}.
Then
$\eta_{\gamma}\restrict |\nabla'|=\eta_{\gamma'}$,
whence in particular  $\eta_{\gamma'}$ is
linear over each simplex of
$|\nabla'|$.
\end{lemma}


  We next prove

\begin{lemma} \label{lemma:independence-bis} Let $\mathcal W=
W_{0},W_{1},\ldots$ be a stellar sequence.  Let $\epsilon_{0}$ be a
$\Delta_{0}$-realization of $W_{0}$ and $\iota_{0}$ be a
$\nabla_{0}$-realization of $W_{0}$.  Then $ \mathcal G(\mathcal W,
\Delta_{0}) \cong \mathcal G(\mathcal W, \nabla_{0}).  $ \end{lemma}

\begin{proof}
Let us write for short $I=\mathcal{I}(\mathcal{W},\Delta_0)$,
$J=\mathcal{I}(\mathcal{W},\nabla_0),$ $P=|\Delta_0|$ and
$Q=|\nabla_0|.$ Then $P\in \mathcal Z(I)$ and $Q\in \mathcal Z(J)$.
In the light of { Theorem \ref{theorem:leo-manu}} it is enough to
exhibit a $\mathbb{Z}$-homeomorphism $\eta$ of $P$ onto $Q$ mapping
$\mathcal{Z}(I)_{\cap P}$ one-one onto $\mathcal{Z}(J)_{\cap Q}.$ By
definition of realization, there is a combinatorial isomorphism $\xi$
of $W_{\Delta_0}$ onto $W_{\nabla_0}$.  By { Lemma
\ref{lemma:unione-jems} (i)}, $\xi$ can be extended to a
$\mathbb{Z}$-homeomorphism $\eta$ of $P$ onto $Q,$ which is linear
over each simplex of $\Delta_0$.  { Lemma \ref{lemma:tre-sei}} now
yields $\mathbb Z$-homeomorphisms \begin{equation}
\label{equation:lorena-i} \eta\,\restrict\, |\Delta_{i}|\, \colon
|\Delta_{i}| \cong_{\mathbb Z} |\nabla_{i}|, \,\,\, i=0,1,2,\ldots,
\end{equation} with $\eta\restrict |\Delta_{i}|$ linear on every
simplex of $\Delta_{i}$.
Since every $X \in \mathcal Z(I)_{\cap P}$ contains some
$|\Delta_i|\in \mathcal O(\mathcal W,\Delta_{0})$ then by
(\ref{equation:lorena-i}), $$ Q\,\,\supseteq \,\,\eta(X)\,\, \supseteq
\,\,\eta(|\Delta_i|) \,=\, |\nabla_{i}| \,\in \,\mathcal O(\mathcal
W,\nabla_{0}), $$ whence by { Proposition \ref{proposition:woj}}, $
\eta(X)\in \mathcal Z(J)_{\cap Q}.  $ Reversing the roles of $\eta$
and $\eta^{-1}$ we conclude that the map $X\mapsto \eta(X)$ sends
$\mathcal Z(I)_{\cap P}$ one-one onto $\mathcal Z(J)_{\cap Q}$, as
desired.  \end{proof}

\subsection*{Proof of Theorem \ref{theorem:surject}}
By { Corollary \ref{corollary:quotient}}
there exists an integer
$n>0$  such that $(G,u)$ is unitally $\ell$-isomorphic to $\McNn/I$
for some $\ell$-ideal $I$ of $\McNn$.
We list the elements of $I$ in a sequence $f_0, f_1, \ldots $ Let
$P_i=\bigcap_{j=0}^i \mathcal Z (f_i),$ for
each $i=0,1,2,\ldots$

Since    $\mathcal Z(f_i)\in \mathcal Z(I)$
and $\mathcal Z(I)$ is closed
under finite  intersections,
   $P_i$ belongs to $\mathcal Z(I).$
Moreover, for each $f\in I$ there is $j = 0,1,2,\ldots$ such that
$P_j\subseteq \mathcal Z(f).$ Thus the filtering family $\{ P_0, P_1,
\ldots\}$ is confluent with $\mathcal Z(I)$, in symbols,
\begin{equation} \label{equation:pyz} \{ P_0, P_1,
\ldots\}\,\,\uparrow\,\, \mathcal Z(I).  \end{equation} By { Lemma
\ref{lemma:support}}, $P_0$ is the support of a regular complex
$\Delta_0.$ Proposition \ref{proposition:extended-dcp} yields a finite
sequence of regular complexes $$ \Delta_{0,0}, \Delta_{0,1} ,\ldots ,
\Delta_{0,k_0} $$  having the following properties:

\begin{enumerate}
\item
    $\Delta_{0,0}=\Delta_0; $
    \item
    for each $t =1,2,\ldots,$
    $\Delta_{0,t}$ is obtained  by  blowing-up
   $\Delta_{0,t-1}$
    at the  Farey mediant of some   $1$-simplex
    $E\in \Delta_{0,t-1}$;
    \item
    $P_1$ is a union of
simplexes of $\Delta_{0,k_0}.$
\end{enumerate}
Let the sequence of regular complexes
$
\Delta_{0,k_0}, \Delta_{0,k_0+1}, \ldots , \Delta_{0,r_0}
$
be obtained by the following procedure:
for each $i>0,$ delete in $\Delta_{0,k_0+i-1}$ a
maximal simplex $T$  which is not contained in $P_1$;
denote by
$\Delta_{0,k_0+i}$  the resulting complex;
stop when no such $T$ exists. Then the sequence of skeletons
$
W_{\Delta_{0,0}}, \ldots,
W_{\Delta_{0,k_0}},\ldots, W_{\Delta_{0,r_0}}
$ is
a finite initial segment of a stellar sequence
and $|\Delta_{0,r_0}|=P_1$.
Let us write   $\Delta_{1,0}$ instead of
$\Delta_{0,r_0}.$
Proceeding inductively, we obtain
a sequence  $\mathcal{S}$
of regular complexes
$$
\mathcal{S}=
\Delta_{0,0},  \ldots , \Delta_{1,0}, \ldots ,
\Delta_{2,0},\ldots,
\Delta_{j,0}, \ldots
$$
such that $P_j = |\Delta_{j,0}|$ for each
$j = 0,1,2,\ldots .$

\medskip To conclude the proof, let $\mathcal W$ be the stellar
sequence given by the skeletons of the regular complexes in $\mathcal
S.$ Let $\rho$ be the trivial $\Delta_{0}$-realization of the skeleton
$W_{\Delta_{0}}$ of $\Delta_{0}$.  Recalling (\ref{equation:pyz}) we
get $$ \mathcal O(\mathcal W, \Delta_{0})= \mathcal O(\mathcal W,
\Delta_{0,0})= \{|\Delta_{0,0}|, \ldots , |\Delta_{1,0}|,\ldots,
\}\uparrow \{P_0, P_1, \ldots\}\uparrow \mathcal Z(I).  $$
Applying {Lemma \ref{lemma:independence-bis} and Proposition
\ref{proposition:tre}} we obtain $$ \mathcal G(\mathcal W) \cong
\mathcal G(\mathcal W, \Delta_{0}) = \McNn/I \cong (G,u), $$ which
concludes the proof.  $\hfill \Box$


\section{Stellar sequences of isomorphic unital $\ell$-groups}
In this section we  give
    simple necessary and sufficient conditions
for two stellar sequences to generate isomorphic unital
$\ell$-groups:

\begin{theorem}
        \label{theorem:equivalence-copy}
For any two stellar
sequences $\mathcal W$ and ${\bar{\mathcal W}}$
let us write
$\mathcal O(\mathcal W)=
|\Delta_{0}|\supseteq |\Delta_{1}|\supseteq \ldots $,
and
${\mathcal O}(\bar{\mathcal W)} =
|\bar{\Delta_{0}}|\supseteq
|\bar{\Delta_{1}}|\supseteq \ldots $
Then the
following conditions are equivalent:
\begin{itemize}
       \item[(i)]
$\mathcal G(\mathcal W) \cong
\mathcal G({\bar{\mathcal W}})$.

\smallskip

\item[(ii)] For some integer $i\geq 0$ there is  a
${\mathbb Z}$-hom\-eo\-mor\-phism $\eta$
of $|\Delta_{i}|$ such that
    $ \{\eta(|\Delta_{i}|),
    \eta(|\Delta_{i+1}|),\ldots\}
\uparrow \mathcal O({\bar{\mathcal W}}).$

    \end{itemize}
    \end{theorem}


    \begin{proof}
Writing  for short $I=
\mathcal I(\mathcal W)$ and $J=\mathcal I({\bar{\mathcal W}})$,
we have
$\mathcal G(\mathcal W)= {\McNm}/{I}$
and
$\mathcal G({\bar{\mathcal W}})= {\McNn}/{J}$.
Thus
\begin{equation}
       \label{equation:luce}
\mathcal O(\mathcal W)
     \subseteq \mathcal Z(I),\,\,
\mbox{ and }\,\, \mathcal O(\mathcal W)
\mbox{ minorizes }\,\,    \mathcal Z(I)
\end{equation}
and
\begin{equation}
       \label{equation:luce2}
\mathcal O(\bar{\mathcal W})
     \subseteq \mathcal Z(J),
\,\,\mbox{ and }\,\, \mathcal O(\bar{\mathcal W})
\mbox{ minorizes }\,\,   \mathcal Z(J).
\end{equation}

\medskip
\noindent
       $(i)\to(ii)$.
{ Theorem \ref{theorem:leo-manu}
(i)$\to$(ii)}  gives    rational
polyhedra  $P\in \mathcal Z(I)$ and
$Q\in \mathcal Z(J)$,
and a
$\mathbb Z$-homeomorphism
    $\eta$ of $P$ onto $Q$ such that
$X\mapsto \eta(X)$
maps
     $\mathcal Z(I)_{\cap P}$
     one-one onto
     $\mathcal Z(J)_{\cap Q}.$
By (\ref{equation:luce}),
there exists an  $i$ such that
$|\Delta_{i}|\subseteq P$.
Define
$D=|\Delta_{i}|$ and
$E=\eta(D).$
Then
    \begin{equation}
        \label{equation:secondo}
     X\mapsto \eta(X)\,\,
     \mbox{ maps }\,\,
     \mathcal Z(I)_{\cap D}\,\,
    \mbox{ one-one onto }\,\,\mathcal Z(J)_{\cap E}.
\end{equation}
In the light of
{ Proposition \ref{proposition:threesets}},
from  $D\in \mathcal O(\mathcal W)$
we get $D \in \mathcal Z(I)_{\subseteq D}
= \mathcal Z(I)_{\cap D}, $
whence  $E\in \mathcal Z(J)_{\subseteq E}
= \mathcal Z(J)_{\cap E} .$
{}Upon writing
$$
\mathcal O(\mathcal W)_{\cap D} =
\{X\cap D\mid X\in \mathcal O(\mathcal W)\}
\mbox{ and }
\mathcal O(\bar{\mathcal W})_{\cap E} =
\{Y\cap E\mid Y\in \mathcal O(\bar{\mathcal W)}\}, $$
  by (\ref{equation:luce2}) we have the confluence
$\mathcal O(\mathcal W)_{\cap D}
    \uparrow \mathcal Z(I)_{\cap D} .$
    Similarly,
    $\mathcal O(\bar{\mathcal W})_{\cap E}
    \uparrow \mathcal Z(J)_{\cap E}$,
whence
$
\eta(\mathcal O(\mathcal W)_{\cap D})\uparrow
\eta(\mathcal Z(I)_{\cap D}).
$
Recalling now  (\ref{equation:secondo}),
we obtain
$\eta(\mathcal Z(I)_{\cap D})=
    \mathcal Z(J)_{\cap E}  \uparrow
    \mathcal O(\bar{\mathcal W})_{\cap E}.
    $
Recalling that  $E\supseteq |\bar{\Delta_{j}}|$, for some
integer
$j\geq 0$ we get
$$ \{\eta(|\Delta_{i}|),
    \eta(|\Delta_{i+1}|),\ldots\}=
    \eta(\mathcal O(\mathcal W)_{\cap D})
\uparrow  \mathcal O(\bar{\mathcal W})_{\cap E}
\uparrow \mathcal O(\bar{\mathcal W}),$$
as required.

%
%
%


    \bigskip
\noindent  $(ii)\to(i)\,\,\,$
Writing
$P =|\Delta_{i}|$ and $Q=\eta(P)$, we have
$
\mathcal Z(I)_{\cap P}   \,\,\uparrow\,\,
\mathcal Z(I) \,\,\uparrow\,\,
\mathcal O(\mathcal W)\,\,\uparrow\,\,
      \mathcal O(\mathcal W)_{\cap P}.
$
By hypothesis,
$
\eta(\mathcal O(\mathcal W)_{\cap P})
\,\,\uparrow\,\, \mathcal O(\bar{\mathcal W})
\,\,\uparrow\,\, \mathcal Z( J).
$
Since $P\in \mathcal O(\mathcal W)$
then $Q\in \mathcal Z(J).$
For every
$Y\in  \mathcal Z( J)_{\cap Q} $ there is
$Z \in   \mathcal O(\mathcal W)_{\cap P}$  such that
$\eta(Z)\subseteq Y.$
    Thus $Z \subseteq
\eta^{-1}(Y) \subseteq P,$ whence
$\eta^{-1}(Y)$ belongs to $\mathcal Z(I)_{\cap P} $.
Reversing the roles of $\eta$ and
$\eta^{-1}$ we see that the
map $X\in \mathcal Z(I)_{\cap P} \mapsto \eta(X)$
is (one-one)  onto $ \mathcal Z(J)_{\cap Q} $.
An application of   {Theorem  \ref{theorem:leo-manu}
(ii)$\to$(i)}
provides the required unital $\ell$-isomorphism
     $\McNm/I\cong \McNn/J $.
\end{proof}
%
%


Our second   characterization depends on the
following

\begin{definition} Two stellar sequences $\mathcal{W}=
W_{0},W_{1},\ldots$ and\ $\mathcal{W}^{\prime}=
W_{0}^{\prime},W_{1}^{\prime},\ldots$ are {\it strongly equivalent} if
there are realizations $\rho:W_{0}\rightarrow\nabla_0$ of $W_{0}$ and
$\rho^{\prime}:W_{0}^{\prime}\rightarrow\nabla_0^{\prime}$ of
$W_{0}^{\prime}$ in the same $d$-cube, whose orbits are confluent, $
\mathcal O (\mathcal W ,\nabla_0) \uparrow \mathcal O (\mathcal
W^{\prime} ,\nabla_0^{\prime}).  $ \end{definition}

\begin{theorem}
\label{theorem:equivalence-bis}
For any
two stellar sequences $\mathcal{W}$ and $\mathcal{W}'$,
$\mathcal G (\mathcal W) \cong \mathcal G (\mathcal
W^{\prime})$  if and only if
there is a stellar sequence
$\mathcal{W}''$ which is strongly equivalent
   to $\mathcal{W}$ and to
$\mathcal{W}'.$
\end{theorem}

\begin{proof}
      ($\Leftarrow$) This is an immediate consequence
      of { Proposition \ref{proposition:tre}} and
{ Lemma \ref{lemma:independence-bis}}.
($\Rightarrow$) Let us write
$\mathcal{W}=W_{0},W_{1},\ldots$ and
$\mathcal{W^{\prime}}
=W'_{0},W'_{1},\ldots.$
Let $\rho_0\colon W_0\rightarrow\Delta_{0}$ and
$\rho'_0\colon W_0^{\prime}\rightarrow\Delta_{0}^{\prime}$
be  the  canonical realizations of $W_{0}$ and
$W_{0}^{\prime}$, respectively in the $m$-  and
in the $n$-cube.
Then
{ Lemma \ref{lemma:meno-star}} yields
commutative diagrams

     \begin{equation}
       \label{equation:infinito1}
\def\normalbaselines{\baselineskip20pt
\lineskip3pt
\lineskiplimit3pt}
\def\mapright#1{\smash{\mathop{\longrightarrow}\limits^{#1}}}
\def\mapdown#1{\Big\downarrow
\rlap{$\vcenter{\hbox{$\scriptstyle#1$}}$}}
\begin{matrix}
W_{0}  & \mapright{\flat_{0}} &
W_{1}  & \mapright{\flat_{1}} &
W_{2}  & \ldots\cr
\mapdown{\rho_{0}} &  &
\mapdown{\rho_{1}} &  &
\mapdown{\rho_{2}} & \cr
     \Delta_{0} & \mapright{\flat_0(\rho_{0})} &
     \Delta_{1} & \mapright{\flat_1(\rho_{1})} &
     \Delta_{2} & \ldots\cr
\end{matrix}
\end{equation}

\noindent and

\begin{equation}
       \label{equation:infinito2}
\def\normalbaselines{\baselineskip20pt
\lineskip3pt
\lineskiplimit3pt}
\def\mapright#1{\smash{\mathop{\longrightarrow}\limits^{#1}}}
\def\mapdown#1{\Big\downarrow
\rlap{$\vcenter{\hbox{$\scriptstyle#1$}}$}}
\begin{matrix}
W'_{0}  & \mapright{\flat'_{0}} &
W'_{1}  & \mapright{\flat'_{1}} &
W'_{2}  & \ldots\cr
\mapdown{\rho'_{0}} &  &
\mapdown{\rho'_{1}} &  &
\mapdown{\rho'_{2}} & \cr
     \Delta'_{0} & \mapright{\flat'_0(\rho'_{0})} &
     \Delta'_{1} & \mapright{\flat'_1(\rho'_{1})} &
     \Delta'_{2} & \ldots\cr
\end{matrix}
\end{equation}
\begin{quote}
\end{quote}
Next,
{ Theorem \ref{theorem:leo-manu}}
yields two polyhedra
$P\in\mathcal{Z(I(W))}$ and
$Q\in\mathcal{Z(I(W^{\prime}))}$, and a
$\mathbb{Z}$-homeomorphism $\eta:P\rightarrow Q$ such
   that the map $X\mapsto \eta(X)$ sends $\mathcal{Z(I(W))}_{\subseteq
P}$ one-one onto $\mathcal{Z(I(W^{\prime}))}_{\subseteq Q}. $
Since $\mathcal O(\mathcal W)$
is confluent with
$\mathcal{Z(I(W))}$
there is $|\Delta_{i^*}|\in\mathcal O (\mathcal W)$ such that
$|\Delta_{i^*}|\subseteq P.$
We will build a sequence
$$
\mathcal{S}=\{\Theta_{0,0},\ldots,\Theta_{0,n_{0}},
\Theta_{1,0},\ldots,\Theta_{1,n_{1}},
\Theta_{2,0},\ldots,\Theta_{2,n_{2}},
\ldots
, \Theta_{j,0}, \ldots, \Theta_{j,n_{j}},
\ldots \}
$$ of regular complexes having the following properties:
\begin{itemize}

       \item[(i)]  for
each $j\geq 0,$
   $\eta$
      {\it respects} $\Theta_{j,n_{j}}$,
      in the sense that $\eta$
      is linear over each simplex
      of $\Theta_{j,n_{j}}$;

    \item[(ii)]
     $\Theta_{j,n_{j}}$ is a subdivision of
$\Delta_{i^*+j},\,\,\,(j=0,1,\ldots)$.
\end{itemize}

\medskip
\noindent{\bf Construction of $\mathcal{S}$.}
Setting  $\Theta_{0,0}=\Delta_{i^*}$,
by { Lemma \ref{lemma:dominio_lineal}}
we have a sequence of  binary Farey blow-ups
$\Theta_{0,0},\Theta_{0,1},\ldots,\Theta_{0,n_{0}}$ such that
$\eta$ respects $\Theta_{0,n_0}$.   Thus
conditions (i) and (ii) are satisfied for
$j=0.$
%
%
Proceeding by induction,  and writing
$k=i^*+j$, let
$\Theta_{j,n_{j}}$ be a subdivision of
$\Delta_{k}$  such that $\eta$
respects $\Theta_{j,n_{j}}$.
The sequence  sequence
$\Theta_{j+1,0},\ldots,\Theta_{j+1,n_{j+1}}$
is now constructed  arguing by cases:

\medskip
\noindent {\it Case 1:}
$\flat_k\colon W_k \rightarrow W_{k+1}$ is a binary
subdivision.
Then  the $\Delta_k$-trans\-for\-m\-a\-tion
$\flat_k(\rho_{k}):\Delta_{k} \rightarrow \Delta_{k+1}$ is a binary
Farey blow-up.  Again,
{ Lemma \ref{lemma:dominio_lineal}} guarantees the
existence of a  sequence
$ \Theta_{j,n_j},\Theta_{j+1,0}, \ldots,\Theta_{j+1,n_{j+1}}$
of binary Farey blow-ups such that
$\Theta_{j+1,n_{j+1}}$ is    a subdivision
of $\Delta_{k+1}.$
   By induction hypothesis,
   $\eta$ respects $\Theta_{j+1,n_{j+1}}$.

\medskip
\noindent {\it Case 2:} $\flat_{k}\colon W_{k}
\rightarrow W_{k+1}$ is a deletion of
a maximal set $M$.
Then by  { Lemma \ref{lemma:meno-star}}  the
$\Delta_k$-transformation
$\flat_k(\rho_{k}):\Delta_{k} \rightarrow
\Delta_{k+1}$ deletes the  corresponding
maximal simplex
$\conv(\rho_{k}(M))$.
Since $\Theta_{j,n_j}$ is
a subdivision of $\Delta_k$,   there is a finite sequence
$\Theta_{j,n_j},\Theta_{j+1,0}, \ldots,\Theta_{j+1,n_{j+1}}$
    such that every complex in the sequence
    is obtained from the
previous one via a deletion of a maximal simplex and
$\Theta_{j+1,n_{j+1}}$ is a subdivision of $\Delta_{k+1}.$
By induction hypothesis,  $\eta$
respects $\Theta_{j+1,n_{j+1}}.$

\medskip
\noindent {\it Case 3:}
$\flat_{k}\colon W_{k} \rightarrow W_{k+1}$ is the identity
map.
Then $\Delta_{k+1}=\Delta_{k}$
and it is sufficient to set $n_{j+1}=0$ and
$\Theta_{j+1,0}=\Theta_{j,n_{j}}.$

\bigskip
\noindent Defining now
$\mathcal{W}''=
W_{\Theta_{0,n_0}},W_{\Theta_{1,0}},\ldots,W_{\Theta
_{1,n_{1}}},\ldots ,$
it follows that    $\mathcal{W}''$
is a stellar sequence.
\bigskip

\noindent
{\it Claim 1:}
$\mathcal W$ is strongly equivalent to $\mathcal W''.$

As a matter of fact,
the trivial realization of the skeleton
of ${\Theta_{0,n_0}}$ yields
the sequence of rational polyhedra
$$
\mathcal O(\mathcal W'',\Theta_{0,n_0})=
\{|\Theta_{0,n_0}|,|\Theta_{1,0}|,\ldots, |\Theta_{1,n_1}|,\ldots,
|\Theta_{i,n_{i}}|,\ldots\}.
$$
By (ii), for all $j=1,2,\ldots,$
$|\Delta_{i^*+j}|=|\Theta_{j,n_{j}}|.$
It follows that
$\mathcal O (\mathcal W)$ is confluent with
$\mathcal O (\mathcal
W'',\Theta_{0,n_0}).$

\bigskip
\noindent
{\it Claim 2:} $\mathcal W'$ is strongly
equivalent to $\mathcal W''.$

Let us write for short $\Theta$ instead of ${\Theta_{0,n_0}}$.
Let $\mathcal{V}=\{v_{1},\ldots v_{k}\}$ be the set of
vertices  of
${\Theta}$.
We define the complex  $\Lambda=\eta(\Theta)$  by the
following stipulation:
$$\conv(\eta(v_{i_{0}}),\ldots\eta(v_{i_{p}}))
\in  \Lambda \Leftrightarrow
   \conv(v_{i_{0}},\ldots v_{i_{p}}) \in
   {\Theta}.$$
Let   $\mathcal{V}'$ be the set
of vertices of  $\Lambda$.
Since  $\eta$ is linear on each simplex of
${\Theta},\,\,$ by { Corollary \ref{corollary:reg_into_reg}}
$\,\,\Lambda$  is a regular complex.
Let the map  $\xi\colon \mathcal{V}\to \mathcal{V'}$
be defined by   $\xi(v)=\eta(v)$.
Since $\eta$ preserves least common denominators
of rational vertices, $\xi$ is a
combinatorial
isomorphism of $W_{{\Theta}}$ onto $W_{\Lambda},$
whence $\xi\colon
W_{\Theta}\rightarrow\Lambda$
is a  $\Lambda$-realization of $W_{\Theta}$,
and  we have
$$\mathcal O(\mathcal W^{\prime\prime},\Lambda)=
\{|\Lambda|, |\Lambda_{1,0}|,\ldots,|\Lambda_{1,n_{1}}|,
\ldots\}=
\{\eta(|\Theta|),  \eta(|\Theta_{1,0}|),
\ldots,\eta(|\Theta_{1,n_{1}}|),\ldots\}.$$
{}From
$|\Delta_{i^*}|=|\Theta|\in \mathcal{Z(I(W))}_{\cap P}$
it follows that  $|\Lambda_{i,j}|$ belongs to
$\mathcal{Z(I(W^{\prime}))}_{\cap Q},$
for all $i,j$.
Therefore,
$\mathcal O(\mathcal W^{\prime})$ minorizes $\mathcal
O(\mathcal W^{\prime\prime},\Lambda).$
Since
$Q\in \mathcal{Z(I(W'))}$ then
$|\Delta_k^{\prime}|\subseteq Q$,
for some $k$.
{}From the hypothesis that $\eta$
maps $\mathcal{Z(I(W))}_{\cap P}$
one-one onto
$\mathcal{Z(I(W^{\prime}))_{\cap Q}},$
it follows  that
$\eta^{-1}\left(  |\Delta_{k+l}^{\prime}|\right)$
belongs to $ \mathcal{Z(I(W))}_{\cap P}$,
for all $l=0,1,2,\ldots, \ $.
In other words, there exists
$j\geq 0$ such that
$\eta^{-1}\left(  |\Delta_{k+l}^{\prime}|\right)  \supseteq|\Delta
_{i^*+j}|=|\Theta_{j,n_{j}}|.$
Summing up,
for each $l\geq 0,$ there is
$j\geq 0$ such
that
$
|\Delta_{k+l}^{\prime}|\supseteq\eta\left(  |\Theta_{j,n_{j}}|\right)
=|\Lambda_{j,n_j} |,
$
which settles our second claim and concludes the proof.
\end{proof}

\section{Examples and Concluding Remarks}

\begin{theorem}
\label{theorem:tavola}
For  $\,\,\mathcal W=W_{0}, W_{1},\ldots$  a stellar sequence,
let  $\mathcal O(\mathcal W)$, $\mathcal I(\mathcal W)$
and
    $\mathcal G(\mathcal W)$ be as given by our main
    construction.
We then have:
\begin{enumerate}
       \item   $\mathcal G(\mathcal W)$
       is finitely presented
       (i.e.,  the $\ell$-ideal
       $\mathcal I(\mathcal W)$ is finitely generated)
       iff $\mathcal W$
       has only finitely many deletion steps.

        \smallskip
    \item $\mathcal G(\mathcal W)$
    is finitely presented and its maximal spectrum
(\cite[10.2.2]{bigkeiwol}, \cite[Section 2]{marmun})
    has dimension
    at most $1$ iff  $\mathcal W$
       has only finitely many deletion steps
       and for all large $i$
    all
    simplexes of $W_{i}$ have at most two elements.

           \smallskip
       \item  $\mathcal G(\mathcal W)$ is
       {\rm simplicial}  \cite[p.191]{effshe}
       (i.e., $\mathcal G(\mathcal W)$ is a finite
       product of unital $\ell$-groups of the form
       $(\mathbb Z\frac{1}{n},1)$)  iff
       for all large $i$  every
       simplex of $W_{i}$, other than
       the empty set, is a  singleton.

       \smallskip
       \item  $\mathcal G(\mathcal W)$
       is archimedean  \cite[2.6.1]{bigkeiwol},   iff every
       rational polyhedron $P$ containing
      the intersection  $\bigcap
    \mathcal O(\mathcal W)$  of the zerosets
      of $\mathcal I(\mathcal W)$
also contains some zeroset   $Q \in \mathcal O(\mathcal W)$.

           \smallskip
       \item  $\mathcal G(\mathcal W)$
       is {\rm local}  (i.e.,   has
      precisely one maximal $\ell$-ideal)
       iff  $\bigcap
    \mathcal O(\mathcal W)$
is a singleton.

           \smallskip
       \item   $\mathcal G(\mathcal W)$
       is unitally $\ell$-isomorphic to a
       unital $\ell$-subgroup of
       $(\mathbb R, 1)$
       iff $\,\,\,\bigcap
    \mathcal O(\mathcal W)$
is a singleton  $\{z\}$
and every rational polyhedron
       $P\ni z$  contains some
       $Q\in
    \mathcal O(\mathcal W)$.

\smallskip
       \item $\mathcal G(\mathcal W)$  is totally ordered
       iff  for some (equivalently, for every)
       $\nabla$-realization of
       $\mathcal W$  (with $\nabla$ a regular complex),
       we have:
       For all rational polyhedra  $P,Q$ with
       $P \cup Q =[0,1]^{n}$  there is  $0\leq i\in \mathbb Z$
    such that either  $|\nabla_{i}|\subseteq
       P$ or $|\nabla_{i}|\subseteq Q$.

    \end{enumerate}
       \end{theorem}

       \begin{proof} We have
\begin{enumerate}
	\item  ($\Leftarrow$)   Let $W_{u}$  be such that
	no step after the $u$th step is a deletion.
All supports $|\Delta_{W_{u+t}}|$ are equal to
	some fixed  rational polyhedron $P$. Let
	$J_{P}$  be the $\ell$-ideal of all functions
	of $\McNn$ that vanish over $|\Delta_{W_{u}}|$.
	By our main construction,
	$\mathcal G(\mathcal W)\cong \McNn/J_{P}
	\cong \McN(P)$, whence
	by \cite[5.2]{marmun},
	$\mathcal G(\mathcal W)$ is finitely presented.

	\medskip
\noindent
($\Rightarrow$) Suppose
	$\mathcal G(\mathcal W)$ is finitely
	presented.  By \cite[5.2]{marmun}
the intersection
$P=\bigcap
	 \mathcal O(\mathcal W)$
	is a rational polyhedron
	and we can identify
    $\mathcal G(\mathcal W)$  and $\McN(P)$.
Further,   $\mathcal G(\mathcal W)$ is archimedean, and
	$\McN(P)
	\cong \McNn/\mathcal I(\mathcal W)$.
	A
    rational polyhedron $Q$ belongs to
    $\mathcal Z(\mathcal I(\mathcal W))$
    iff $Q\supseteq P$ iff $Q$ is minorized
    by some $R\in \mathcal O(\mathcal W)$.
    In particular, taking
    $Q=P$  we see that
    $P$ belongs to $\mathcal O(\mathcal W)$,
    say  $P=|\Delta_{i}|$. It follows that
    $|\Delta_{i}|=|\Delta_{i+1}|=\ldots$.

	\smallskip
\item  For some  rational polyhedron $P$,
$\mathcal G(\mathcal W)$
can be identified with  $\McN(P)$.
Our assumption about the
maximal spectrum of
$\mathcal G(\mathcal W)$  is equivalent to
saying that for all large $i$
the $i$th element of $\mathcal O(\mathcal W)$
is the same rational polyhedron $P$, and
$\dim P\leq 1$.  This means that  no
simplex of $W_{i}$  has three elements.

\smallskip
	\item  This is a particular case of (2).

\smallskip
	\item  By \cite[2.3]{marmun}
the    unital
	$\ell$-group $\mathcal G(\mathcal W) \cong
	\McNn/\mathcal I(\mathcal W)$ is archimedean iff
	it is unitally $\ell$-isomorphic to
	$\McN(X)$ for some closed
	set $X\subseteq [0,1]^{n}.$
Identifying $\mathcal G(\mathcal W)$ and $\McN(X)$,
	the $\ell$-ideal
	$\mathcal I(\mathcal W)$ turns out to
	coincide with the
	set of functions in  $\McNn$  that vanish over
	$\bigcap \mathcal O(\mathcal W)$, and
$X=
\bigcap \mathcal O(\mathcal W)
=
\bigcap \{\mathcal Z(f) \mid  f\in \mathcal I(\mathcal W)\}.$
Thus if $\mathcal G(\mathcal W)$ is archimedean
every rational polyhedron $P$ containing
$\bigcap \mathcal O(\mathcal W)$  automatically
contains some $Q\in \mathcal O(\mathcal W)$.
If $\mathcal G(\mathcal W)$ is not archimedean then
some $f\in \McNn$  vanishing over $\bigcap \mathcal O(\mathcal W)$
does not belong to $\mathcal I(\mathcal W)$.
Thus the rational polyhedron $\mathcal Z(f)$
contains $\bigcap \mathcal O(\mathcal W)$ but
    does not contain any $Q\in \mathcal O(\mathcal W)$.

		\smallskip
	\item  Trivial.

	\smallskip
	\item  From the time-honored
	H{\"o}lder theorem,
	\cite[2.6]{bigkeiwol},
	it follows that
	a unital
	$\ell$-group
is unitally $\ell$-isomorphic to
$(\mathbb R,1)$ iff  its only
$\ell$-ideal is $\{0\}$, iff it is
    archimedean and local.

	\smallskip
	\item
Write $J$ instead of
$\mathcal I(\mathcal W)$.
As is well known,
\cite[2.4.1]{bigkeiwol},
	the  quotient
	$\mathcal G(\mathcal W)=\McNn/J$
	is totally ordered iff
	$f \wedge g = 0$ implies
	$f \in J \,\,\mbox{or}\,\,
       g\in J$.   Since
       (by { Proposition \ref{proposition:woj}})
       ideals are uniquely
       determined by their zerosets,
       this is equivalent to saying that
       for any two disjoint
rational polyhedra  $P= \mathcal Z(f)$
and $Q=\mathcal Z(g)$ such that
$  P \cup Q = [0,1]^{n}$,  either
$P$ or $Q$ is the zeroset of
some function $f\in J$.
In other words, either
$P$ or $Q$
is minorized by some
       $R \in \mathcal O(\mathcal W)$.
       The same proof holds for any
       $\nabla$-realization.

\end{enumerate}
\end{proof}

\begin{examples}
       \label{examples:table}
We give a few instances of the map
$\mathcal W\mapsto \mathcal G(\mathcal W)$:

\begin{enumerate}
    \item    $\mathcal W =$     the constant sequence
       $W,W,W,\ldots$ where
        $W=(\mathcal V,\Sigma,\omega)$,
$\mathcal V=\{0,1\},$  $\Sigma={\rm powerset}(\mathcal V),$
and $\omega=1.$
       Then
       $\mathcal G(\mathcal W)\cong \McN([0,1])$.



           \smallskip
       \item     $\mathcal W = $   the
       constant sequence $W,W,\ldots$
       with $W=(\mathcal V,\Sigma,\omega)$
       such that $\mathcal V=\{w_{1},\ldots,w_{k}\}$,
       $\Sigma=\{\emptyset, \{w_{1}\},\ldots,\{w_{k}\}\}$
       and $\omega(w_{i})=n_{i}, \,\,\,
       i=1,\ldots,k.$   Then
       $\mathcal G(\mathcal W)$ is the simplicial
       group
$
       (\mathbb Z\frac{1}{n_{1}},1)\times\cdots
\times(\mathbb Z\frac{1}{n_{k}},1).$

           \smallskip
       \item
        $\mathcal W =$   the
       constant sequence $W,W,\ldots$, where
       $W$ is the skeleton of a
    regular complex     $\Delta$
       whose support is contained in $[0,1]^{n}.$
       Then   $\mathcal G(\mathcal W)\cong\McN(|\Delta|)$.

           \smallskip
       \item    Let  $W_{0}=(\mathcal V,\Sigma,\omega)$
    where  $\mathcal V=\{0,1\}$,  $\Sigma=
    {\rm powerset}(\mathcal V)$,
       and $\omega(v)=n$  for all $v\in \mathcal V. $
       Define $W_{i}$ inductively by
       $$
       W_{i}=
    \begin{cases}
       \mbox{subdivide the only existing 2-element set in } W_{i-1}
       &\mbox{ if } i \equiv 1 \mod 3\cr
           \mbox{delete the only 2-element set $X$ of
	$W_{i-1}  $ with $0\notin X$ }
	&\mbox{ if } i\equiv 2 \mod 3\cr
\mbox{delete the only maximal singleton $\{s\}$ of $W_{i-1}$}
	&\mbox{ if } i \equiv 0 \mod 3.\cr
       \end{cases}
       $$
Let $\mathcal W=W_{0},W_{1},\ldots$.
Then
$\mathcal G(\mathcal W)\cong
(\mathbb Z \otimes_{\rm lex}
       \mathbb Z, (n,0))$  is the unital
$\ell$-group    given by the free abelian
group $\mathbb Z^{2}$  with
lexicographic ordering  and with
the element $(n,0)$
as the order-unit.

               \smallskip
       \item  For  $\xi$  an irrational in
       the unit real interval
       define the sequence of regular complexes
	$\Delta_{0},
	\Delta_{1}, \Delta_{2},  \ldots$
	as follows: let
	$\Delta_{0}$  the complex given  by the unit interval $[0,1]$
	 and its faces;
	 and for $i>0$,
	    $$
\Delta_{i}=
    \begin{cases}
       \mbox{do a Farey blow-up
       of the only  1-simplex in } \Delta_{i-1}
       &\mbox{ if } i \equiv 1 \mod 3\cr
           \mbox{delete  the only 1-simplex $X$ of
	$\Delta_{i-1}  $ with $\xi \notin X$ }
	&\mbox{ if } i \equiv 2 \mod 3\cr
\mbox{delete the only maximal 0-simplex $\{s\}$ of
$\Delta_{i-1}$}
	&\mbox{ if } i \equiv 0 \mod 3.\cr
       \end{cases}
       $$
Let $W_{0}=(\mathcal V, \Sigma, \omega)$ be defined by
     $\mathcal V=\{0,1\}$,  $\,\,\Sigma =
{\rm powerset}(\mathcal V)$,
       and $\omega=1,$ over all of $\mathcal V$.
       Observe that $W_{0}$ is the skeleton of
       $\Delta_{0}$.
       Let  $W_{i}$  be similarly obtained as the
       skeleton of
       $\Delta_{i}$.
Then $\mathcal G(\mathcal W)
\cong
(\mathbb Z\xi+\mathbb Z,1)$  is the
unital $\ell$-group generated in $\mathbb R$
by $\xi$, with $1$ as the order-unit.

           \end{enumerate}
\end{examples}


\subsection*{Final Remarks} Up to isomorphism, every stellar sequence
$\mathcal W$ determines a unique AF C*-algebra $A=A_{\mathcal W}$ via
the map $$ \mathcal W\mapsto \mathcal G(\mathcal W) \mapsto K_{0}^{-1}
(\mathcal G(\mathcal W)), $$ where $K_{0}(A)$ is the unital dimension
group of $A$, \cite{goo-af}.  Combining Elliott classification
\cite{ell, goo-af} with Theorem {\ref{theorem:surject}, we get that
(up to isomorphism) the range of the map $\mathcal W\mapsto
A_{\mathcal W}$ coincides with the class of unital AF C*-algebras $A$
whose dimension group $K_{0}(A)$ is lattice-ordered and finitely
generated.  Various important AF C*-algebras existing in the
literature belong to this class, including the Behnke-Leptin algebra
with a two-point dual \cite{behlep}, (whose $K_{0}$-group is as in
Example (4) above) the Effros-Shen algebras \cite{effshe} (whose
$K_{0}$-groups are those given in Example (5) above), and various
algebras considered in \cite{cigellmun}, the universal AF C*algebra
$\mathfrak M_{1}$ of \cite{mun-adv} (= the algebra $\mathfrak A$ of
\cite{boc}) (whose $K_{0}$-group is $\McN([0,1])$).  Theorem
\ref{theorem:equivalence-copy} provides a simple criterion to
recognize when two stellar sequences $\mathcal W$ and $\mathcal W'$
determine isomorphic AF C*-algebras $A_{\mathcal W}$ and $A_{\mathcal
W'}$.  This criterion is a simplification of the equivalence criterion
for Bratteli diagrams, \cite[2.7]{bra}.



\begin{thebibliography}{10}
\bibitem{ale} J.W. Alexander,
The combinatorial theory of complexes,
{\it Annals of Mathematics,} 31:292-320, 1930.

\bibitem{behlep}
H. Behncke, H. Leptin,
C*-algebras with a two-point dual,
{\it Journal of
Functional  Analysis},   10:330-335, 1972.


\bibitem{bigkeiwol} {A.
     Bigard, K. Keimel, and S. Wolfenstein,}
{\it Grou\-pes et Anneaux
     R\'{e}ticul\'{e}s},   Lecture Notes in
Mathematics, Springer-Verlag, Berlin, volume 608, 1971.


 \bibitem{boc}
{F. Boca},
An AF algebra associated with the Farey tessellation,
{\it  Canad. J. Math},
{60}:975--1000, 2008.


\bibitem{bra}
O. Bratteli,  Inductive limits of
finite-dimensional C*-algebras,
{\it Trans. Amer. Math.
Soc.,}   171:195-234, 1972.






\bibitem{cigellmun}
R. Cignoli, G.A. Elliott, D.Mundici,
Reconstructing $C^{*}$-algebras from
their Murray von Neumann orders,
{\it Advances in Mathematics,} 101:166-179, 1993.






\bibitem{effshe}
E.G. Effros, C.L. Shen, Approximately finite
C*-algebras and continued fractions,
{\it Indiana  J. Math.},   29:191-204, 1980.

\bibitem{ell}
G.A. Elliott, On the classification of
inductive limits of sequences of semisimple
finite-dimensional algebras,
{\it J. Algebra},   38:29-44, 1976.




\bibitem{ewa}
G. Ewald,  {\it Combinatorial convexity and
algebraic geometry},
Springer-Verlag, New York, 1996.





\bibitem{glahol}
A.M.W. Glass, W.C. Holland, (Eds.)
{\it Lattice-ordered groups,}
Kluwer Academic Publishers, 1989.


 %
 %



\bibitem{goo-af}
{K.R. Goodearl},
{\it Notes on real and complex C$^{*}$-algebras},
Birk{\"a}user, Boston, Inc. Shiva
Math. Series, 5,  1982.





\bibitem{grulek}
{P.M. Gruber,  C.G. Lekkerkerker},
{\it Geometry of numbers}, 2nd edition,
North-Holland, Amsterdam,
1987.


 %
 %
 %


    \bibitem{marmun}
      V. Marra,  D. Mundici,
      The Lebesgue state of a unital
      abelian lattice-ordered group,
      {\it Journal of Group Theory,}
10:655-684, 2007.





\bibitem{mun-jfa}
{ D. Mundici,}
Interpretation of {A}{F}
${C}\sp \ast$-algebras
in {{\L}}ukasiewicz sentential calculus,
    {\it Journal of Functional
Analysis,} 65:15-63, 1986.



\bibitem{mun-adv}
{ D. Mundici},   Farey stellar subdivisions,
ultrasimplicial groups, and
$K_{0}$ of AF $C^{*}$-al\-ge\-bras,
{\it Advances in Mathematics,} 68:23-39,
1988.

    %




    %





    %
    %




     \bibitem{mun-dcds}
     { D. Mundici,}
     The Haar theorem for lattice-ordered abelian
   groups with order-unit,
    {\it Discrete and Continuous
     Dynamical Systems,}  21:537-549,  2008.





    %
    %






 %
 %
 %
 %




\bibitem{sta}
J.R. Stallings,
{\it Lectures on Polyhedral Topology,}
Tata Institute of Fundamental Research,
Mumbay, 1967.



\bibitem{wlo}
{J. W{\l}odarczyk},
Decompositions of birational toric maps
in blow-ups and blow-downs.
{\it Trans. Amer. Math. Soc.,}
349:373-411, 1997.


\end{thebibliography}
\end{document}